\documentclass[11pt]{article}
\usepackage{amssymb,amsmath}
\newtheorem{thm}{Theorem}
\newtheorem{cor}[thm]{Corollary}
\newtheorem{prop}[thm]{Proposition}
\newtheorem{lemma}[thm]{Lemma}
\newtheorem{definition}[thm]{Definition}

\newcommand {\bD} {\mathbb {D}}
\newcommand {\bC} {\mathbb {C}}

\newcommand {\bZ} {\mathbb {Z}}

\newcommand {\cD} {\mathcal {D}}
\newcommand {\cK} {\mathcal {K}}

\newcommand {\cF} {\mathcal {F}}
\newcommand {\fT} {\mathfrak T}

\newcommand {\pa} {\partial}
\newcommand {\Tr} {\text {Tr}}
\newcommand {\Ker} {\text {Ker}}
\newcommand {\Cker} {\text {Coker}}
\newcommand {\Ran} {\text {Ran}}
\newcommand {\Ind} {\text {Index}}
\newcommand {\Dim} {\text {dim}}

\def\mapright#1{\smash{\mathop{\longrightarrow}\limits^{#1}}}
\newcommand {\proof} {\noindent{\it Proof. }}

\begin{document}

\title{Dirac operators on noncommutative manifolds with boundary.}

\author{\ Alan L. Carey \\ Mathematical Sciences Institute \\
Australian National University \\ Canberra ACT0200 
AUSTRALIA\ 
\\
\\ S\l awomir Klimek \\ Department of Mathematics \\
Indiana University Purdue University Indianapolis\\ 402 N. Blackford St.
 Indianapolis, IN 46202
\\
\\
Krzysztof P. Wojciechowski \footnote{Unfortunately K. P. Wojciechowski died in June 2008 while this paper was being prepared and his ideas were not fully incorporated. We, his friends and collaborators, sorely miss his inspiration.} 
\\ Department of Mathematics \\
Indiana University Purdue University Indianapolis\\ 402 N. Blackford St. Indianapolis, IN 46202}

\date{December 30, 2008}
\maketitle

\begin{abstract}
\vskip 0.3cm
\noindent We study an example of an index problem for a Dirac-like operator subject to
 Atiyah-Patodi-Singer  boundary conditions on a noncommutative manifold with boundary, namely the quantum unit disk.

\end{abstract}

\section{Introduction}\label{intesec}

The objective of this paper is to construct simple noncommutative analogs of Dirac operators on manifolds with boundary and to study noncommutative versions of Atiyah-Patodi-Singer (APS) boundary conditions. The example we describe is essentially an analog of the d-bar operator on the two dimensional unit disk.

The noncommutative (quantum) disk we consider was described in \cite{KL1}. It is defined as the $\bC^*$-algebra generated by certain weighted unilateral shifts. It has natural analogs of angular and radial coordinates so that calculations can be made quite explicit.

For the theory of APS and other global elliptic boundary conditions for Dirac type operators on manifolds with boundary see \cite{BW}, and references therein. 

A related study of an example of APS boundary conditions in the context of noncommutative geometry is contained in \cite{CPR}. There, APS boundary conditions are considered for a `noncommutative cylinder'. One starts with non-commutative algebras $A,B$, and assumes there is an unbounded $( {A, B})$ Kasparov module $(X,\cal D)$
where $\cal D$ is a regular operator on the module $X$ having discrete spectrum. This is `boosted' to an unbounded Kasparov module $(L^2([0,\infty))\otimes X, \hat{\cal D})$ 
 with APS boundary conditions (from the spectral splitting of $\cal D$) for a certain regular operator $\hat{\cal D}$.  Thus in \cite{CPR} one has a noncommutative boundary algebra and the cylinder does not introduce additional noncommutativity. The situation in the present paper is different in that the `noncommutative algebra with boundary' is already given. It has  a commutative boundary algebra and the noncommutativity arises from the `interior'. The problem we resolve here is also different: namely we find a naturally defined unbounded analogue of the d-bar operator which is Fredholm for the APS boundary conditions determined by the spectral splitting of a classical  Dirac operator on the boundary.

The paper is organized as follows. In Section \ref{classicsec} we present the details of
the index calculations for d-bar like operators in weighted $L^2$- spaces on
the two-dimensional disk. Section \ref{disksec} contains the definition and our study of the properties of noncommutative analogs of the operators from the previous section on the quantum disk. Finally in Section \ref{examplesec} we look in some detail at one example
inspired by references \cite{KL1} and \cite{BKLR}.

\section{APS}\label{classicsec}

In this section we review the basic elements of the elliptic theory for the APS like boundary conditions of simple Dirac type operators on the disk in two dimensions.
In particular we will explicitly compute the polar decomposition of the operator, its kernel and then look at an appropriate integration by parts formula. This is followed by explicit index calculations.

First we introduce some notation. Our example is on the unit disk:
\begin{eqnarray*}
\bD:=\{z\in\bC: |z|\leq 1\},
\end{eqnarray*}
the boundary of which is
\begin{eqnarray*}
\partial \bD=\{z\in\bC: |z|= 1\}\cong S^1.
\end{eqnarray*}
We will use polar coordinates $z=\rho e^{i\phi}$ and in particular the Fourier transform in variable $\phi$.

We want to study the following operator, a multiple of the d-bar operator:
\begin{eqnarray*}
D=F(\rho)\frac{\pa}{\pa \bar z}
\end{eqnarray*}
with the coefficient $F(\rho)>0$ and $F(1)=2$. The coefficient $F(\rho)$ also
appears in the weight of the Hilbert space, and in fact can be easily entirely eliminated from the discussion. The main reason we keep it is that its noncommutative analog in the next section is essential and its presence here makes the formulas in both sections look very similar.

In polar coordinates the operator is:
\begin{eqnarray}
D=\frac{F(\rho)}{2\rho}e^{i\phi}\left(\rho\frac{\pa}{\pa \rho}+i\frac{\pa}{\pa \phi}\right)\label{poldec1}
\end{eqnarray}
Let us mention that this decomposition shows that the operator $D$ has a slightly different structure than the one assumed in \cite{APS} and \cite{BW}. (The difference is essentially by the factor $e^{i\phi}$.) This clearly has no effect on the
rest of the theory.

We need more notation: the restriction to the boundary map is
\begin{eqnarray*}
r:C(\bD)\to C(\pa\bD),\ \ \ rf(\phi)=f(1\cdot e^{i\phi}).
\end{eqnarray*}
We have the exact sequence:
\begin{eqnarray}
0\ \mapright{}\ C_0(\bD)
\ \mapright{}\ C(\bD)\ \mapright{r}\ C(\partial \bD)\ \mapright{}\ 0,
\label{exseq1}
\end{eqnarray} 

Choose any ``extension" map $e:C(\pa\bD)\to C(\bD)$, such that $ef(\rho e^{i\phi})=f(\phi)$ near $\rho=1$.
Clearly $r\circ e=1$.

Restriction of the operator $D$ is
\begin{eqnarray}
rDe=ie^{i\phi}\frac{\pa}{\pa \phi}\label{Dboundary}
\end{eqnarray}
Again, in the usual APS setup one considers operators which are $i\frac{\pa}{\pa \phi}$ on the boundary,
and then uses spectral properties of that operator to set up boundary conditions on $D$. We do the same here,
ignoring the $e^{i\phi}$ factor.

For future comparisons we write the formula for $D$ in Fourier transform.
For $f(z)=\sum\limits_{n\in\bZ}f_n(\rho)e^{in\phi}$ we have
\begin{eqnarray}
Df(z)=\frac{F(\rho)}{2\rho}\sum_{n\in\bZ}(\rho f_n'(\rho)-nf_n(\rho))
e^{i(n+1)\phi}\label{DFT1}
\end{eqnarray}

Clearly, the kernel of $D$ consists of holomorphic functions on $\bD$. Consequently
\begin{eqnarray}
r(\Ker D)=\text{span}\left(e^{in\phi}\right)_{n\geq 0}\label{kernel1}
\end{eqnarray}

Next we look at $D$ in a weighted $L^2$ space and in particular we compute the formal adjoint of $D$, and the corresponding integration by parts formula.

Consider the following $L^2$ inner product on functions on $\bD$ with the weight equal to the reciprocal of $F(\rho)$:
\begin{eqnarray}
(f,g)=\int\bar f(z) g(z) F(\rho)^{-1}\frac{dz\wedge d\bar z}{-2i\pi}\label{innprod1}
\end{eqnarray}

The corresponding Hilbert space is denoted by $L^2_F(\bD)$, and let 
$H^1_F(\bD)\cong H^1(\bD)$ be the first Sobolev space on $\bD$. It is well known that the restriction map $r$ extends to $H^1(\bD)$.

We want to find the adjoint of $D$ in $L^2_F(\bD)$. To this end we
compute as follows, using the Stokes' theorem and for $f,g\in L^2_F(\bD)$:
\begin{eqnarray*}
(Df,g)&=&\int_\bD\overline{F(\rho)\frac{\pa f}{\pa \bar z}}\, g(z) F(\rho)^{-1}\frac{dz\wedge d\bar z}{-2i\pi}=\\
&=&\int_\bD d(\bar f(z)g(z))\wedge\frac{d\bar z}{-2i\pi}-\int_\bD\bar f(z)
F(\rho)\frac{\pa g(z)}{\pa  z}F(\rho)^{-1}\frac{dz\wedge d\bar z}{-2i\pi}=\\
&=&\int_{\pa\bD}\bar f(z)g(z)\frac{d\bar z}{-2i\pi}+\int_\bD\bar f(z)
\bar Dg(z)F(\rho)^{-1}\frac{dz\wedge d\bar z}{-2i\pi}=\\
&=&\int_0^{2\pi}\overline{rf(\phi)}\,rg(\phi)e^{-i\phi}\frac{d\phi}{2\pi}+(f,\bar Dg).
\end{eqnarray*}
where
\begin{eqnarray*}
\bar D=-F(\rho)\frac{\pa}{\pa  z}
\end{eqnarray*}

Notice that the kernel of the $\bar D$ operator consists of anti-holomorphic functions and so:
\begin{eqnarray}
r(\Ker \bar D)=\text{span}\left(e^{in\phi}\right)_{n\leq 0}\label{kernel2}
\end{eqnarray}

We now impose non-local APS like boundary conditions on $D$ in order to make it Fredholm. For this we need the following
notation. Let $P_N$ be the projection in $L^2(S^1)$ onto $\text{span}\left(e^{in\phi}\right)_{n\leq N}$. The point is that $P_N$ are spectral projections of the (modified) boundary operator $i\frac{\pa}{\pa \phi}$, see (\ref{Dboundary}). 
Define $D_{P_N}$ to be the operator $D$ on the domain
\begin{eqnarray*}
\text{dom}\,(D_{P_N}):=\{ f\in H^1(\bD):rf\in\Ran\, P_N \}.
\end{eqnarray*}

The following is a routine index theorem.
\begin{thm}\label{indexthm}
{\ }

The operator $D_{P_N}$ is an unbounded Fredholm operator. Its adjoint 
$D_{P_N}^*$ is the operator $\bar D$ with the domain
\begin{eqnarray*}
\text{dom}\,(D_{P_N}^*):=\{ f\in H^1(\bD):e^{-i\phi}rf\in\Ker\, P_N \}.
\end{eqnarray*}
Moreover:
\begin{eqnarray*}
\Ind\, D_{P_N}=N+1.
\end{eqnarray*}

\end{thm}
\proof 
The theorem follows from the general theory - see \cite{APS} and \cite{BW}. The index calculation is as follows.
\begin{eqnarray*}
\Dim\,\Ker(D_{P_N})&=&\#\{n:\ 0\leq n\leq N\}=\\
&=&\begin{cases} 
0& \text{if $N<0$}\\
N+1& \text{if $N\geq 0$.}
\end{cases}
\end{eqnarray*}
Similarly:
\begin{eqnarray*}
\Dim\,\Ker(\bar D_{P_N})&=&\#\{n:\ N+1< n\leq 0\}=\\
&=&\begin{cases} 
-(N+1)& \text{if $N<0$}\\
0& \text{if $N\geq 0$,}
\end{cases}
\end{eqnarray*}
and the formula for the index follows. $\square$

\section{NC disk}\label{disksec}

In this section we describe an analog of the constructions and results of the previous section when
the usual unit disk $\bD$ is replaced by the noncommutative disk. We will recycle some notation from the previous
section: for example symbols $r,e,D,\bar D$ will denote similar but different objects than before.

According to \cite{KL1}, the Toeplitz algebra $\fT$ plays the role of the noncommutative disk. The Toeplitz algebra $\fT$ is defined, in terms of generators and relations, as the universal unital $\bC^*$-algebra generated by a single generator $U$ subject to relation $U^*U=1$.
Alternatively, and more concretely, it is the unital $\bC^*$-algebra generated by the unilateral shift $U:l^2\to l^2$, $Ue_k:=e_{k+1}$, where $\{e_k\}_{k\geq 0}$ is the canonical basis in $l^2$, see \cite{F}.
The formula for $U^*$ is
\begin{eqnarray*}
U^*e_k=\begin{cases} 
0& \text{if $k=0$}\\
e_{k-1}& \text{if $k\geq 1$}
\end{cases}\label{basis}
\end{eqnarray*}

We will extensively use the following ``label" operator in $l^2$:
\begin{eqnarray*}
Ke_k:=ke_k.
\end{eqnarray*}
In particular, the commutation relation for functions of $K$ with $U$ will be frequently used:
\begin{eqnarray}
U^*f(K)=f(K+1)U^*.\label{comrel1}
\end{eqnarray}

We have the exact sequence, similar to (\ref{exseq1}):
\begin{eqnarray}
0\ \mapright{}\ \cK
\ \mapright{}\ \fT\ \mapright{r}\ C(S^1)\ \mapright{}\ 0,
\label{exseq2}
\end{eqnarray}
where the restriction map $r$ will be described later. Here $\cK$ is the algebra of compact operators. Consequently
we think of $\fT$ as a noncommutative manifold with boundary, where the boundary is $S^1$, same as in the commutative case.

There is a natural extension map $e:C(S^1)\to\fT$ that can be conveniently described in terms of Toeplitz operators as follows. We identify $l^2$ with the subspace of $L^2(S^1)$ spanned by $e^{in\phi}$, $n\geq 0$ and let $\Pi$ be the orthogonal projection onto it. For $f\in C(S^1)$, let $M_f$ be the multiplication by $f$ operator in $L^2(S^1)$. Then:
\begin{eqnarray*}
e(f):=\Pi M_f\in\fT.
\end{eqnarray*}

To make things as explicit as possible as well as to define ``generalized" functions we will work with Fourier series expansion. Leaving aside any convergence questions,
if $a\in\fT$ then we have the following expansion:
\begin{eqnarray}
a=\sum_{n\geq 1}f_n(K)(U^*)^n + \sum_{n\geq 0}U^ng_n(K)\label{FSeries},
\end{eqnarray}
where
\begin{eqnarray}
f_n, g_n\in\cF:=\{f:\bZ_{\geq 0}\to \bC: f^b:=\lim_{k\to\infty}f(k)\ \text{exists}\}.
\label{bvalue}
\end{eqnarray}

The above restriction map $r$ is, in Fourier transform, given by:
\begin{eqnarray*}
r\left(\sum_{n\geq 1}f_n(K)(U^*)^n + \sum_{n\geq 0}U^ng_n(K)\right)=
\sum_{n\geq 1}f_n^be^{-in\phi} + \sum_{n\geq 0}g_n^be^{in\phi}.
\end{eqnarray*}

Similarly, the extension map can be expressed in Fourier transform as:
\begin{eqnarray*}
e\left(\sum_{n\geq 1}f_ne^{-in\phi} + \sum_{n\geq 1}g_ne^{in\phi}\right)=
\sum_{n\geq 1}f_n(U^*)^n + \sum_{n\geq 0}U^ng_n.
\end{eqnarray*}

More generally we say that $a$ is a distribution if it is given by a formal power series (\ref{FSeries})
where coefficients $f_n,g_n$ are arbitrary arithmetic functions. 
Most of the formulas below make sense for such distributions. 

We will start considering convergence issues only when the appropriate Hilbert space is introduced and that Hilbert space convergence is the only one needed in this paper. 

The basic idea now is to replace derivatives with commutators and so to consider operators of the form $a\mapsto P[Q,a]$, where $P,Q$ are possibly unbounded operators affiliated with $\fT$.
We will make simple choices to achieve the following objectives that would  make it possible to impose APS like boundary conditions, prove the Fredholm property and compute the index. The analog of the operator $D$ of the previous section should:
\begin{itemize}
\item have a polar decomposition like (\ref{poldec1})
\item reduce to $ie^{i\phi}\frac{\pa}{\pa \phi}$ on the boundary, as in (\ref{Dboundary})
\item have infinite dimensional kernel, the restriction of which to the boundary should coincide with (\ref{kernel1})
\item have an integration by parts formula similar to the one before.
\end{itemize}

With that in mind, we make the following definition.

\begin{definition}\label{Ddef}
{\ }

Let $a$ be a distribution in the sense described above. We set
\begin{eqnarray}
D(a):=A(K)[UB(K),a]
\end{eqnarray}
with the following assumptions on functions $A(k),B(k)$:
\begin{enumerate}
\item $A(k)>0$, $\sum\limits_{k\geq 0}A(k)^{-1}<\infty$
\item $B(k)> 0$, $k\to B(k)$ is increasing, $\lim\limits_{k\to\infty}B(k)=1$
\item $\lim\limits_{k\to\infty}A(k)((B(k+1)-B(k))=1$, 
$\lim\limits_{k\to\infty}\frac{A(k+1)}{A(k)}=1$
\end{enumerate}
\end{definition}

In particular $A(K)$ is unbounded with a trace class resolvent, while $B(K)\in\fT$. The inverse of $A(K)$ will be used as the weight for defining the Hilbert space.  We use the first two conditions for all the analysis of $D$ including the index computation. The last condition is needed only to recover  $ie^{i\phi}\frac{\pa}{\pa \phi}$ on the boundary.

\begin{lemma}\label{limitlemma}
{\ }

For $n\geq 0$ we have
\begin{eqnarray*}
&\lim\limits_{k\to\infty}&A(k+1)((B(k)-B(k+n))=n\\
&\lim\limits_{k\to\infty}&A(k)((B(k)-B(k+n))=-n,
\end{eqnarray*}
an similarly,
\begin{eqnarray*}
&\lim\limits_{k\to\infty}&A(k+n+1)((B(k+n)-B(k))=-n\\
&\lim\limits_{k\to\infty}&A(k+n)((B(k+n)-B(k))=n.
\end{eqnarray*}
\end{lemma}
\proof  
This is a consequence of a simple telescoping argument and the condition $3.$ of the above definition.
$\square$

A simple example of functions satisfying all the requirements of the Definition \ref{Ddef} is
described in Section \ref{examplesec}.

We start the analysis of  $D$ by considering its form in Fourier transform, obtained after using the commutation relation (\ref{comrel1}):
\begin{eqnarray*}
D(a)&=&\sum_{n\geq 1}A(K)UB(K)f_n(K)(U^*)^n-
\sum_{n\geq 1}A(K)f_n(K)(U^*)^nUB(K)+\\
&+&\sum_{n\geq 0}A(K)UB(K)U^ng_n(K)-
\sum_{n\geq 0}A(K)U^ng_n(K)UB(K)=\\
&=&\sum_{n\geq 1}A(K)\left(B(K-1)f_n(K-1)-B(K+n-1)f_n(K)\right)(U^*)^{n-1}\\
&+&\sum_{n\geq 0}U^{n+1}A(K+n+1)\left(B(K+n)g_n(K)-B(K)g_n(K+1)\right)
\end{eqnarray*}
In the formula above we set $B(-1)=f_n(-1)=0$ to make sense of $B(K-1)$ and 
$f_n(K-1)$. Additionally $U(U^*)^n=\chi(K)(U^*)^{n-1}$ where
\begin{eqnarray*}
\chi(k)=\begin{cases} 
0& \text{if $k=0$}\\
1& \text{if $k\geq 1$}
\end{cases}\label{chidef}
\end{eqnarray*}
However since $f_n(K-1)\chi(K)=f_n(K-1)$, $\chi(K)$ was left out in the formula for $D$.

An analog of a polar decomposition (\ref{poldec1}) and (\ref{DFT1}) is
\begin{eqnarray*}
D(a)&=&
\sum_{n\geq 1}A(K)B(K+n-1)\left(f_n(K-1)-f_n(K)\right)(U^*)^{n-1}+\\
&+&\sum_{n\geq 0}U^{n+1}A(K+n+1)B(K)\left(g_n(K)-g_n(K+1)\right)+\\
&+&\sum_{n\geq 1}A(K)\left((B(K-1)-B(K+n-1)\right)f_n(K-1)(U^*)^{n-1}+\\
&+&\sum_{n\geq 0}U^{n+1}A(K+n+1)\left(B(K+n)-B(K)\right)g_n(K)=\\
&=&D_{radial}(a)+D_{angular}(a)
\end{eqnarray*}
Intuitively, the operator $U$ plays the role of an angular coordinate while $K$ is the (discrete) radial coordinate with the difference operator replacing the corresponding derivative.

\begin{prop}
{\ }

We have, just as in (\ref{Dboundary})
\begin{eqnarray*}
rDe=ie^{i\phi}\frac{\pa}{\pa \phi}
\end{eqnarray*}
\end{prop}
\proof  
For $f=\sum\limits_{n\geq 1}f_n\,e^{-in\phi} + \sum\limits_{n\geq 0}g_n\,e^{in\phi}$ the radial part of $Def$ vanishes. 
What is left is
\begin{eqnarray*}
Def&=&\sum_{n\geq 1}A(K)\left((B(K-1)-B(K+n-1)\right)\chi(K)f_n(U^*)^{n-1}+\\
&+&\sum_{n\geq 0}U^{n+1}A(K+n+1)\left(B(K+n)-B(K)\right)g_n
\end{eqnarray*}
Consequently:
\begin{eqnarray*}
rDef&=&e^{i\phi}\sum_{n\geq 1}\lim_{k\to\infty}A(k)((B(k-1)-B(k+n-1))
f_ne^{-in\phi}+\\
&+&e^{i\phi}\sum_{n\geq 0}\lim_{k\to\infty}A(k+n+1)(B(k+n)-B(k))
g_ne^{in\phi}=\\
&=&ie^{i\phi}\frac{\pa f}{\pa \phi}
\end{eqnarray*}
because the limits inside are equal to $\mp n$ by Lemma \ref{limitlemma}.
$\square$

Now we compute the kernel of D.

\begin{prop}\label{kernelprop}
{\ }

If a distribution $a=\sum\limits_{n\geq 0}f_n(K)(U^*)^n + \sum\limits_{n\geq 1}U^ng_n(K)$ is in the kernel of $D$, then there are constants $a_n$ such that
\begin{eqnarray*}
a=\sum_{n\geq 0} a_n(UB(K))^n.
\end{eqnarray*}
\end{prop}
\proof  
Firstly, $A(K)$ is invertible so $Da=0$ implies $[UB(K),a]=0$.
In Fourier transform this means
\begin{eqnarray*}
[UB(K),f_n(K)(U^*)^n]=0
\end{eqnarray*}
for $n\geq 1$ and
\begin{eqnarray*}
[UB(K),U^ng_n(K)]=0
\end{eqnarray*}
for $n\geq 0$.

The first of those equations means
\begin{eqnarray*}
B(K-1)f_n(K-1)-B(K+n-1)f_n(K)=0
\end{eqnarray*}
This can be solved recursively as follows:
\begin{eqnarray*}
0-B(n-1)f_n(0)&=&0\\
B(0)f_n(0)-B(n)f_n(0)&=&0\\
B(1)f_n(1)-B(n+1)f_n(1)&=&0\\
\ldots&=&0
\end{eqnarray*}
It is clear that the only solution is $f_n(k)=0$.

Now we will analyze the second equation.
Proceeding as before we get
\begin{eqnarray*}
B(K+n)g_n(K)-B(K)g_n(K+1)=0
\end{eqnarray*}
Explicitly, this means
\begin{eqnarray*}
B(n)g_n(0)-B(0)g_n(1)&=&0\\
B(n+1)g_n(1)-B(1)g_n(2)&=&0\\
\ldots&=&0
\end{eqnarray*}
Solving this recursively, we see that
\begin{eqnarray*}
g_n(k)&=&g_n(0)\frac{B(n)B(n+1)\ldots B(n+k-1)}
{B(0)B(1)\ldots B(k-1)}=\\
&=&g_n(0)\frac{B(k)B(k+1)\ldots B(n+k-1)}
{B(0)B(1)\ldots B(n-1)}=\\
&=&\text{const}\,B(k)B(k+1)\ldots B(n+k-1),
\end{eqnarray*}
where the const does not depend on $k$.
Comparing the above with
\begin{eqnarray*}
(UB(K))^n=U^nB(K)B(K+1)\ldots B(K+n-1)
\end{eqnarray*}
proves the proposition.
$\square$

As a corollary we get
\begin{eqnarray*}
r(\Ker D)=\text{span}\left(e^{in\phi}\right)_{n\geq 0},
\end{eqnarray*}
just as in (\ref{kernel1}).

Now we consider the adjoint of $D$. As will be demonstrated below, it is given by the formula:
\begin{eqnarray}
\bar D(a):=A(K)[B(K)U^*,a]
\end{eqnarray}

The analysis of $\bar D$ is completely analogous to that of $D$. Below we work out the details of the
polar decomposition of $\bar D$, its restriction to the boundary, and its kernel.

It is convenient to rearrange the Fourier expansion
\begin{eqnarray}
a=\sum_{n\geq 0}f_n(K)(U^*)^n + \sum_{n\geq 1}U^ng_n(K),
\end{eqnarray}
so that, compared to (\ref{FSeries}), $f_0(K)=g_0(K)$. 
With this convention, the formula for $\bar D$ in Fourier transform is:
\begin{eqnarray*}
&\,&\bar D(a)=\sum_{n\geq 0}A(K)B(K)U^*f_n(K)(U^*)^n-
\sum_{n\geq 0}A(K)f_n(K)(U^*)^nB(K)U^*\\
&\,&+\sum_{n\geq 1}A(K)B(K)U^*U^ng_n(K)-
\sum_{n\geq 1}A(K)U^ng_n(K)B(K)U^*=\\
&\,&=\sum_{n\geq 0}A(K)\left(B(K)f_n(K+1)-B(K+n)f_n(K)\right)(U^*)^{n+1}+\\
&\,&+\sum_{n\geq 1}U^{n-1}A(K+n-1)\left(B(K+n-1)g_n(K)-B(K-1)g_n(K-1)\right)
\end{eqnarray*}
This leads to the following polar decomposition of $\bar D$:
\begin{eqnarray*}
&\bar D(a)&=
\sum_{n\geq 0}A(K)B(K)\left(f_n(K+1)-f_n(K)\right)(U^*)^{n+1}+\\
&+&\sum_{n\geq 1}U^{n-1}A(K+n-1)B(K+n-1)\left(g_n(K)-g_n(K-1)\right)+\\
&+&\sum_{n\geq 0}A(K)\left((B(K)-B(K+n)\right)f_n(K)(U^*)^{n+1}+\\
&+&\sum_{n\geq 1}U^{n-1}A(K+n-1)\left(B(K+n-1)-B(K-1)\right)g_n(K-1)=\\
&=&\bar D_{radial}(a)+\bar D_{angular}(a).
\end{eqnarray*}

Computation of the restriction of $\bar D$ to the boundary can now be done as follows.

\begin{prop}
{\ }

We have
\begin{eqnarray}
r\bar De= ie^{-i\phi}\frac{\pa}{\pa \phi}.
\end{eqnarray}
\end{prop}
\proof  
For $f=\sum\limits_{n\geq 1}f_ne^{-in\phi} + \sum\limits_{n\geq 0}g_ne^{in\phi}$ the radial part of $\bar Def$ vanishes as before. Consequently:
\begin{eqnarray*}
\bar Def&=&\sum_{n\geq 0}A(K)\left((B(K)-B(K+n)\right)f_n(U^*)^{n+1}+\\
&+&\sum_{n\geq 1}U^{n-1}A(K+n-1)\left(B(K+n-1)-B(K-1)\right)\chi(K)g_n,
\end{eqnarray*}
and
\begin{eqnarray*}
r\bar Def&=&e^{-i\phi}\sum_{n\geq 0}\lim_{k\to\infty}A(k)((B(k)-B(k+n))
f_ne^{-in\phi}+\\
&+&e^{-i\phi}\sum_{n\geq 1}\lim_{k\to\infty}A(k+n-1)(B(k+n-1)-B(k-1))
g_ne^{in\phi}=\\
&=&ie^{-i\phi}\frac{\pa f}{\pa \phi}
\end{eqnarray*}
because the limits inside are equal to $\mp n$ by Lemma \ref{limitlemma}.
$\square$

Now we compute the kernel of $\bar D$.
\begin{prop}
{\ }

If a distribution $a=\sum\limits_{n\geq 0}f_n(K)(U^*)^n + \sum\limits_{n\geq 1}U^ng_n(K)$ is in the kernel of $\bar D$, then there are constants $a_n$ such that
\begin{eqnarray*}
a=\sum_{n\geq 0} a_n(B(K)U^*)^n.
\end{eqnarray*}
\end{prop}
\proof  
It is enough to study the equation $[B(K)U^*,a]=0$.
Taking the adjoint and using Proposition \ref{kernelprop} proves the statement.
$\square$

As a corollary we get
\begin{eqnarray*}
r(\Ker \bar D)=\text{span}\left(e^{in\phi}\right)_{n\leq 0}
\end{eqnarray*}
just as in (\ref{kernel2}).

Now we are going to consider $D$ and $\bar D$ in a Hilbert space that is an analog of the weighted $L^2$ space of the previous section.
For $a,b\in\fT$ we define (compare with (\ref{innprod1})):
\begin{eqnarray*}
(a,b)_A:=\Tr\left(A(K)^{-1} ba^*\right),
\end{eqnarray*}
and denote by $L_A^2(\fT)$ the completion of $\fT$ with respect to the above scalar product. We use the notation $||a||$ for the norm of $a$ in $\fT$, and $||a||_A$ for the norm of $a$ in $L_A^2(\fT)$.
 
Alternatively, using Fourier transform, the Hilbert space $L_A^2(\fT)$ can be defined as the space of
distributions $a$ with finite norm $||a||_A$ given by;
\begin{eqnarray*}
||a||_A^2&=&\sum_{n\geq 1}\Tr\left(\overline {f_n(K)}f_n(K)A(K)^{-1}\right)+\\
&+&\sum_{n\geq 0}\Tr\left(\overline {g_n(K)}g_n(K)A(K+n)^{-1}\right)
\end{eqnarray*}

Because both $D$ and $\bar D$ make sense as operators on distributions, we will simply consider them in $L_A^2(\fT)$ on the maximal domains:
\begin{eqnarray*}
\text{Dom}_{max}(D):=\{a\in L_A^2(\fT):\ Da \in L_A^2(\fT)\},
\end{eqnarray*}
and similarly for $\bar D$.

As usual the key to the properties of an unbounded Fredholm operator is provided by its parametrix. For example, a parametrix for $D$ can be obtained by solving $Da=b$ for $a$. 
This is easily done in Fourier transform as follows. Writing
$a=\sum\limits_{n\geq 1}f_n(K)(U^*)^n 
+ \sum\limits_{n\geq 0}U^ng_n(K)$ and 
$b=\sum\limits_{n\geq 0}p_n(K)(U^*)^n + \sum\limits_{n\geq 1}U^nq_n(K)$, 
the equation $Da=b$ becomes:
\begin{eqnarray}
B(K-1)f_n(K-1)-B(K+n-1)f_n(K)=A(K)^{-1}p_{n-1}(K)
\label{paraeq1}
\end{eqnarray}
and
\begin{eqnarray}
B(K+n)g_n(K)-B(K)g_n(K+1)=A(K+n-1)^{-1}q_{n+1}(K)
\label{paraeq2}
\end{eqnarray}
In components the first equation (\ref{paraeq1}) reads:
\begin{eqnarray*}
0-B(n-1)f_n(0)&=&A(0)^{-1}p_{n-1}(0)\\
B(0)f_n(0)-B(n)f_n(0)&=&A(1)^{-1}p_{n-1}(1)\\
B(1)f_n(1)-B(n+1)f_n(1)&=&A(2)^{-1}p_{n-1}(2)\\
\ldots&=&\ldots
\end{eqnarray*}
Solving this recursively, we see that
\begin{eqnarray*}
f_n(k)&=&-\frac{1}{B(k+n-1)}\,\frac{p_{n-1}(k)}{A(k)}-\\
&-&\frac{B(k-1)}{B(k+n-1)B(k+n-2)}\,\frac{p_{n-1}(k-1)}{A(k-1)}-\\
&-&\ldots - \frac{B(k-1)B(k-2)\ldots B(0)}{B(n+k-1)B(n+k-2)\ldots B(n-1)}
\,\frac{p_{n-1}(0)}{A(0)}=\\
&=&-\sum_{j=0}^k\frac{B(j)B(j+1)\ldots B(j+n-2)}{B(k)B(k+1)\ldots B(k+n-1)}
\,\frac{p_{n-1}(j)}{A(j)}
\end{eqnarray*}
where we arranged all the ``B" terms to have common denominator.

Now we will analyze the equation (\ref{paraeq2}).
Proceeding as before we get
\begin{eqnarray*}
B(n)g_n(0)-B(0)g_n(1)&=&A(n+1)^{-1}q_{n+1}(0)\\
B(n+1)g_n(1)-B(1)g_n(2)&=&A(n+2)^{-1}q_{n+1}(1)\\
\ldots&=&\ldots
\end{eqnarray*}
Solving this recursively, with $g_n(0)$ to be chosen later, we compute that
\begin{eqnarray*}
g_n(k)&=& - \sum_{j=0}^{k-1} \frac{B(k)B(k+1)\ldots B(k+n-1)}{B(j)B(j+1)\ldots B(j+n)}
\,\frac{q_{n+1}(j)}{A(n+1+j)}+\\
&+& g_n(0)\frac{B(k)B(k+1)\ldots B(n+k-1)}
{B(0)B(1)\ldots B(n-1)}.
\end{eqnarray*}
Anticipating the APS - like boundary conditions we make the choice of $g_n(0)$
such that $\lim\limits_{k\to\infty}g_n(k)=0$. Using $\lim\limits_{k\to\infty}B(k)=1$
yields the following formula for $g_n(0)$:
\begin{eqnarray*}
g_n(0)=\sum_{j=0}^\infty \frac{B(0)B(1)\ldots B(n-1)}{B(j)B(j+1)\ldots B(j+n)}
\,\frac{q_{n+1}(j)}{A(n+1+j)},
\end{eqnarray*}
so that $g_n(k)$ becomes
\begin{eqnarray*}
g_n(k)=\sum_{j=k}^\infty \frac{B(k)B(k+1)\ldots B(k+n-1)}{B(j)B(j+1)\ldots B(j+n)}
\,\frac{q_{n+1}(j)}{A(n+1+j)}.
\end{eqnarray*}
The convergence issues are handled exactly as in the Proposition \ref{Qprop} below.

For convenience we just copy the formulas above to make the following definition.

\begin{definition}
{\ }

With the above notation for $a,b$ we set:
\begin{eqnarray*}
Qb:=\sum_{n\geq 1}(Qp)_n(K)(U^*)^n + \sum_{n\geq 0}U^n(Qq)_n(K)
\end{eqnarray*}
where
\begin{eqnarray*}
(Qp)_n(k):=-\sum_{j=0}^k\frac{B(j)B(j+1)\ldots B(j+n-2)}{B(k)B(k+1)\ldots B(k+n-1)}
\,\frac{p_{n-1}(j)}{A(j)}
\end{eqnarray*}
and
\begin{eqnarray*}
(Qq)_n(k):=\sum_{j=k}^\infty \frac{B(k)B(k+1)\ldots B(k+n-1)}{B(j)B(j+1)\ldots B(j+n)}
\,\frac{q_{n+1}(j)}{A(n+1+j)}
\end{eqnarray*}
\end{definition}

We have the following properties of the operator $Q$:
\begin{prop}\label{Qprop}
{\ }

\begin{enumerate}
\item For every $a\in L_A^2(\fT)$, $DQa=a$.
\item $Q$ is a bounded operator in $L_A^2(\fT)$.
\end{enumerate}
\end{prop}
\proof Notice that the calculations leading to the definition of $Q$ give exactly $DQa=a$.
Alternatively, the proof of part 1 can also be obtained by a direct calculation.

To show that $Q$ is bounded we first use the assumption that $B(k)$ is positive and increasing so that
$\frac{B(j)}{B(k)}\leq 1$ for $j\leq k$, and $\frac{1}{B(k)}\leq\frac{1}{B(0)}$.
It follows that
\begin{eqnarray*}
|(Qp)_n(k)|&\leq&\frac{1}{B(0)}\sum_{j=0}^k\frac{|p_{n-1}(j)|}{A(j)}\leq\\
&\leq&\frac{1}{B(0)}\left(\sum_{j=0}^\infty\frac{1}{A(j)}\right)^{1/2}
\left(\sum_{j=0}^\infty\frac{|p_{n-1}(j)|^2}{A(j)}\right)^{1/2}\\
\end{eqnarray*}
and similarly for $|(Qq)_n(k)|$.
This implies the following estimate on the norm:
\begin{eqnarray*}
||Qb||_A\leq \frac{1}{B(0)}\left(\sum_{j=0}^\infty\frac{1}{A(j)}\right)||b||_A,
\end{eqnarray*}
and the proposition is proved. $\square$

As a corollary we get:
\begin{cor}\label{bvcorr}
(existence of boundary value) 

If $a\in\text{Dom}_{max}(D)$ is given by the Fourier series (\ref{FSeries}), then
the following limits exist for every $n$: $f_n^b:=\lim\limits_{k\to\infty}f_n(k)$ and 
$g_n^b:=\lim\limits_{k\to\infty}g_n(k)$.
\end{cor}
\proof  
If $a\in\text{Dom}_{max}(D)$ then the above Proposition \ref{Qprop} 
implies that $QDa-a$ is in the kernel of $D$. This, together with the kernel calculation of Proposition \ref{kernelprop}, means that there is (a unique) $b\in  L_A^2(\fT)$ and constants $a_n$ such that
\begin{eqnarray*}
a=Qb+\sum_{n\geq 0} a_n(UB(K))^n.
\end{eqnarray*}
It follows that
\begin{eqnarray*}
f_n^b=\lim\limits_{k\to\infty}f_n(k)=
-\sum_{j=0}^\infty B(j)B(j+1)\ldots B(j+n-2)\,\frac{p_{n-1}(j)}{A(j)},
\end{eqnarray*}
where the convergence issues are handled exactly as in the Proposition \ref{Qprop} above.
Our choice of $Q$ implies that:
\begin{eqnarray*}
g_n^b=\lim\limits_{k\to\infty}g_n(k)=a_n,
\end{eqnarray*}
and so both types of limits exist.
$\square$

Consider now the problem of finding a parametrix for $\bar D$. To do that we try to solve for $a$ the equation $\bar D(a)=A(K)[B(K)U^*,a]=b$. After conjugation and rearrangement of terms, the equation is equivalent to $D(a^*)=-A(K)b^*A(K)^{-1}$, which we have already solved leading to the definition of $Q$.

Using $a^*=\sum\limits_{n\geq 1}\bar g_n(K)(U^*)^n
+ \sum\limits_{n\geq 0}U^n\bar f_n(K)$ and 
$b^*=\sum\limits_{n\geq 1}\bar q_n(K)(U^*)^n + \sum\limits_{n\geq 0}U^n\bar p_n(K)$, 
and the formulas for the components of $Q$ we get:
\begin{eqnarray*}
\bar Qb:=\sum_{n\geq 0}(\bar Qp)_n(K)(U^*)^n + \sum_{n\geq 1}U^n(\bar Qq)_n(K)
\end{eqnarray*}
where
\begin{eqnarray*}
(\bar Qp)_n(k):=-\sum_{j=k}^\infty \frac{B(k)B(k+1)\ldots B(k+n-1)}{B(j)B(j+1)\ldots B(j+n)}
\,\frac{p_{n+1}(j)}{A(j)}
\end{eqnarray*}
and
\begin{eqnarray*}
(\bar Qq)_n(k):=\sum_{j=0}^k\frac{B(j)B(j+1)\ldots B(j+n-2)}{B(k)B(k+1)\ldots B(k+n-1)}
\,\frac{q_{n-1}(j)}{A(n-1+j)}
\end{eqnarray*}

It is easy to establish the following properties of the operator $\bar Q$:
\begin{prop}\label{Qbarprop}
{\ }

\begin{enumerate}
\item For every $a\in L_A^2(\fT)$, $\bar D\bar Qa=a$.
\item $\bar Q$ is a bounded operator in $L_A^2(\fT)$.
\item If $a\in\text{Dom}_{max}(\bar D)$, then for every $n$: $f_n^b:=\lim\limits_{k\to\infty}f_n(k)$ and 
$g_n^b:=\lim\limits_{k\to\infty}g_n(k)$ exist.
\end{enumerate}
\end{prop}

The next step is to derive the integration by parts formula for $D$. 
\begin{prop}
{\ }

For $a\in \text{Dom}_{max}(D)$ and $b\in \text{Dom}_{max}(\bar D)$ we have
\begin{eqnarray}
(Da,b)_A=(a,\bar Db)_A-\int_0^{2\pi}\overline{r(a)(\phi)}\,r(b)(\phi)e^{-i\phi}
\frac{d\phi}{2\pi}
\end{eqnarray}
Here both $r(a) = \sum_{n\geq 1}f_n^be^{-in\phi} + \sum_{n\geq 0}g_n^be^{in\phi}$
and  $r(b)$, given by a similar formula, exist by  Corollary \ref{bvcorr} and Proposition \ref{Qbarprop}.
\end{prop}
\proof  
We have the following formal calculation:
\begin{eqnarray*}
(Da,b)_A&=&\Tr\left(A(K)^{-1} b\,(Da)^*\right)=\Tr\left(b\,\{UB(K)a-aUB(K)\}^*\right)=\\
&=&\Tr\left(\{B(K)U^*b-bB(K)U^*\}\,a^*\right)=\Tr\left(A(K)^{-1}\bar Db\,a^*\right)=\\
&=&(a,\bar Db)_A\\
\end{eqnarray*}
The problem with this calculation is that we used the cyclic property of the trace to an operator which is not trace class. To proceed we write 
$\Tr(a)=\lim\limits_{N\to\infty}\sum\limits_{k=0}^N(e_k,ae_k)$ in the above formulas.
While we skip the details of this somewhat lenghtly calculation, we mention that to get the boundary terms we use the Abel's lemma for sequences:
\begin{eqnarray*}
\sum_{k=0}^n a_k(b_{k+1}-b_k)=a_{n+1}b_{n+1}-a_0b_0+
\sum_{k=0}^n b_{k+1}(a_{k+1}-a_k)
\end{eqnarray*}
This lemma applied to diagonal operators, i.e. functions of $K$, gives the following formula, obtained after taking the limits, whenever they exist:
\begin{eqnarray*}
&&\Tr((f(K-1)-f(K))g(K))=\\
&=&\Tr(f(K)(g(K+1)-g(K)))-
\left(\lim_{k\to\infty}f(k)\right)\,\left(\lim_{k\to\infty}g(k)\right).
\end{eqnarray*}
But the limits do exist because of Corollary \ref{bvcorr} and Proposition \ref{Qbarprop}.
All other details are straightforward.
$\square$

We are now in position to state the main result of this paper.
We impose non-local APS like boundary conditions on $D$ in order to make it Fredholm. For this we again recycle the notation of the previous section: $P_N$ is the projection in $L^2(S^1)$ onto $\text{span}\left(e^{in\phi}\right)_{n\leq N}$. We define $D_{P_N}$ to be the operator $D$ on the domain
\begin{eqnarray*}
\text{dom}(D_{P_N})=\{ a\in \text{Dom}_{max}(D):r(a)\in\Ran\, P_N \}
\end{eqnarray*}

With the above notation we have:
\begin{thm}\label{mainthm}
{\ }

The operator $D_{P_N}$ is an unbounded Fredholm operator. Its adjoint 
$D_{P_N}^*$ is the operator $\bar D$ with the domain
\begin{eqnarray*}
\text{dom}( D_{P_N}^*):=\{ a\in \text{Dom}_{max}(\bar D)  :e^{-i\phi}r(a)\in\Ker\, P_N \}. 
\end{eqnarray*}
Moreover:
\begin{eqnarray}
\Ind\, D_{P_N}=N+1.
\end{eqnarray}

\end{thm}
\proof 
We have already done all the hard work.
The first part of the  theorem follows from the fact that we can construct the parametrix
$Q_{P_N}$ such that $D_{P_N}Q_{P_N}=1-P_{\Cker\, D_{P_N}}$ and
$Q_{P_N}D_{P_N}=1-P_{\Ker\, D_{P_N}}$, where $P_{\Ker\, D_{P_N}}$ and $P_{\Cker\, D_{P_N}}$ are orthogonal projections onto corresponding subspaces.
This is done in a routine way by modifying the operator $Q$ on finite dimensional subspaces: namely the kernel and the cokernel of $D_{P_N}$. 

The evaluation of the adjoint of $D_{P_N}$ and its domain is a consequence of the integration by parts formula.

Finally, the calculation of the index is exactly the same as in the commutative case: see the proof of Theorem \ref{indexthm}.
$\square$

\section{Example}\label{examplesec}

In this section we review some aspects of the quantum unit disk  $C(\bD_\mu)$ of \cite{KL1}. Calculus on the disk gives a natural example of the operator $D$ of the previous section.

For $0 < \mu < 1$ we let $C(\bD_\mu)$ denote the universal unital $\bC^*$-algebra generated by two 
elements $z$ and $\bar z$ which are adjoint to each other, with the following relation:
\begin{eqnarray*}
[\bar z, z] =\mu(I-z\bar z)(I-\bar zz).	
\end{eqnarray*}
Such an universal $\bC^*$-algebra is defined in the following way. If $a$ is a polynomial in $z,\bar z$ we define its norm
as the supremum of $||\pi(a)||$ over all Hilbert space representations $\pi$
satisfying the relation. One verifies that this defines a  sub-$\bC^*$-norm and the corresponding
completion mod the null space gives the universal $\bC^*$-algebra.

Recall that the group of biholomorphisms $SU(1, 1)/\bZ_2$ of the unit disk
$\bD = \{z\in\bC: |z| \leq 1\}$ consists of fractional transformations
\begin{eqnarray}
z\to (az+b)(\bar bz+\bar a)^{-1},\hskip 3mm  |a|^2-|b|^2=1.	\label{action}
\end{eqnarray}
The point of the above definition of the quantum unit disk is that the mapping (\ref{action})
defines also an action of $SU(1, 1)$ on the quantum disk $C(\bD_\mu)$. 

We have the following structure theorem:

\begin{thm}\label{structhm}(see \cite{KL1})

Let, as before, $\{e_k\}$ be the canonical basis in $l_2$, $k=0,1,2\ldots$, and let
$z:l_2\to l_2$ be the following weighted unilateral shift: 

\begin{eqnarray*}
ze_k= \left\{ \frac{(k+1)\mu}{1+(k+1)\mu}\right\}^{1/2}\,e_{k+1},\ \   
n\geq 0	
\end{eqnarray*}
with the adjoint equal to
\begin{eqnarray*}
\bar ze_k = \begin{cases} 0&  k=0,\cr
\left\{ \frac{k\mu}{1+k\mu}\right\}^{1/2}\,e_{k-1}, &k\geq 1.\cr	
\end{cases}
\end{eqnarray*}

Then
$C(\bD_q)\cong \bC^*(z,\bar z)$, where $\bC^*(z,\bar z)$ is the unital $\bC^*$-algebra generated by $z,\bar z$.

\end{thm}

The following is a simple consequence of the structure theorem.
\begin{cor}
{\ }

The $\bC^*$-algebras $C(\bD_\mu)$ are isomorphic to each other and
for every $\mu$, $0<\mu< 1$, and we have
$C(\bD_\mu)\cong \fT$, where $\fT$ is the  Toeplitz algebra. 
\end{cor}

Notice that the commutator $\bar z z-z\bar z =\mu(I-z\bar z)(I-\bar zz)$
is an invertible operator (with an unbounded inverse). Indeed we have
\begin{eqnarray*}
(\bar z z-z\bar z)e_k=\frac{\mu}{(1+k\mu)(1+(k+1)\mu)}\,e_k.
\end{eqnarray*}

Define formally (unbounded) operators $\cD$, $\bar \cD$ on $C(\bD_\mu)$ by
\begin{eqnarray*}
\cD a:=(\bar z z-z\bar z)^{-1}[a,z]
\end{eqnarray*}
\begin{eqnarray*}
\bar \cD a:=(\bar z z-z\bar z)^{-1}[\bar z,a]
\end{eqnarray*}
In a sense the operators $\cD$, $\bar \cD$ are the analogs of the usual complex partial derivatives on the quantum unit disk $C(\bD_\mu)$ because we have the familiar relations: 
\begin{eqnarray}
\cD(1)=0,\ \cD(z)=0,\ \cD(\bar z)=1,\label{Drel1}
\end{eqnarray}
\begin{eqnarray}
\bar \cD(1)=0,\ \bar \cD(z)=1,\ \bar \cD(\bar z)=0.\label{Drel2}
\end{eqnarray}

Let
\begin{eqnarray*}
A(k):=\frac{\mu}{(1+k\mu)(1+(k+1)\mu)}
\end{eqnarray*}
and
\begin{eqnarray*}
B(k)=\left\{ \frac{(k+1)\mu}{1+(k+1)\mu}\right\}^{1/2}
\end{eqnarray*}
so that $(\bar z z-z\bar z)^{-1}=A(K)$ and $z=UB(K)$.
It is a simple exercise to verify that:
\begin{enumerate}
\item $A(k)>0$, $\sum\limits_{k\geq 0}A(k)^{-1}<\infty$
\item $B(k)> 0$, $k\to B(k)$ is increasing, $\lim\limits_{k\to\infty}B(k)=1$
\item $\lim\limits_{k\to\infty}A(k)((B(k+1)-B(k))=1$, 
$\lim\limits_{k\to\infty}\frac{A(k+1)}{A(k)}=1$
\end{enumerate}
This verifies that $D:=-\cD$ and $\bar D:=\bar \cD$ are examples of the operators of the previous section.

Another very similar example of a pair of operators $\cD$ and $\bar\cD$ satisfying
relations (\ref{Drel1}) and (\ref{Drel2}) was constructed in \cite{K} for a different, so called
$q$-deformation of the disk related to quantum groups. The difference between the two examples is that the third condition on coefficients $A(k)$, $B(k)$ is not satisfied for the $q$-disk.
Most of the results of this paper still apply for those operators.

\end{document}